\documentclass[12pt,letterpaper]{article}
\usepackage{amsmath}
\usepackage{amssymb}
\usepackage{enumitem}
\usepackage{mathtools}
\usepackage{amsthm}
\usepackage[utf8]{inputenc}
\usepackage[mathscr]{euscript}
\usepackage{pst-node}
\usepackage{stackengine}
\usepackage{fullpage}
\usepackage[semicolon,authoryear]{natbib}
\usepackage[colorlinks=true,linkcolor=blue,citecolor=blue]{hyperref}
\usepackage{booktabs,tabularx}

\usepackage{tikz-cd}
\usepackage{bm}
\usepackage{geometry}
\usepackage[T1]{fontenc}
\usepackage{titling}
\usepackage{threeparttable}

\usepackage{tcolorbox}
\usepackage{listings}
\usepackage{longtable}

\usepackage{xargs}
\usepackage[colorinlistoftodos,prependcaption,textsize=tiny]{todonotes}
\newcommandx{\sgcmt}[2][1=]{\todo[linecolor=red,backgroundcolor=red!25,bordercolor=red,#1]{#2}}

\definecolor{DarkGreen}{RGB}{30,180,30}

\newif\ifbackmatter%
\newcommand{\backmatter}{\global\backmattertrue}%

\usepackage[title]{appendix}%

\theoremstyle{definition}

\theoremstyle{remark}

\newtheorem*{remark}{Remark}

\tikzcdset{
arrow style=tikz,
diagrams={>={Straight Barb[scale=0.8]}}
}

\DeclareMathOperator\Div{Div}
\DeclareMathOperator\Curl{Curl}
\DeclareMathOperator\cdiv{div}

\newcommand{\interior}[1]{%
  {\kern0pt#1}^{\mathrm{o}}%
}

\title{A time-domain adversarial example demonstrating the failure of $H^1$ elements in $a$-$\phi$ electromagnetic simulation}
\author{Timothy Carlson \& Sanjay Govindjee\\
University of California, Berkeley\\
\texttt{timothy.carlson@berkeley.edu, s\_{}g@berkeley.edu}}

\begin{document}

\maketitle

\begin{abstract}

Frequently, $H^1$ continuous finite elements are used in computational electromagnetics despite it being widely known that such an interpolation is not guaranteed to furnish reliable approximations of solutions to general initial/boundary value problems. It is mathematically more natural to use $H(\Curl)$ elements (edge elements)  than $H^1$ elements, since the former form part of a discrete complex that respects the differential-geometric structure of Maxwell’s equations while the latter do not. However, most work discussing the failure of $H^1$ elements for electromagnetics is esoteric and concrete examples demonstrating the clear advantage of edge elements are generally restricted to quasi-static, frequency-domain, or eigenvalue problems. In this work, we summarize essential facts regarding $H^1$ and $H(\Curl)$ elements and present a physically feasible yet intentionally adversarial time-domain eddy-current example demonstrating the failure of gauge-penalized $H^1$ elements in an $a$-$\phi$ formulation compared to explicit gauge-free mixed scalar/edge $H^1/H(\Curl)$ elements. The results indicate that edge elements can be used to furnish reasonable results for general loading conditions, whereas penalized $H^1$ elements fail (e.g., converge to an incorrect state) even with mesh refinement. 
    
\end{abstract}


\section{Introduction}\label{a:c}

Although the advantages and useful properties of $H(\Curl)$ edge elements are well known and have been thoroughly studied \citep[see, e.g.,][]{Zaglmayr:06}, many works continue to use $H^1$ vector-valued elements in time-domain, $a$-$\phi$ electromagnetics applications (see, e.g., \cite{MORENOMATEOS:23}, \cite{Stewart:26}, \cite{Dorn:21}, and \cite{LUCARINI:22}). It is not unreasonable to assume that this continued prevalence may owe itself to the complexity of the analytic explanations of the failure of $H^1$ elements compared to $H(\Curl)$ edge elements and further to the benign nature of many benchmark problems, which fail to clearly elucidate the underlying issues. Additionally, while there exist numerous examples demonstrating that incorrect eigenvalues are obtained with $H^1$ elements (see, e.g., \cite{Boffi:99} and \cite{Cendes:91}), the time domain eddy current examples present in the literature seem to focus on the comparative accuracy of global responses (see, e.g., \cite{Kuczmann:08}, \cite{kamireddy:22}, and \cite{Nakata:90}). Even when edge elements are deemed superior in such cases, the examples do not elucidate the fact that structural issues plaguing $H^1$ elements can lead to obviously spurious solutions in the time domain as well. Finally, many texts seemingly frame the choice of vector function spaces as an opportunity for \emph{optional} improvement rather than as a necessary means of avoiding issues caused by an inherent analytical deficiency. Given this setting, it is understandable why many, especially those who are new to this particular class of problems or are approaching the field from a multiphysics perspective, may incorrectly view an $H^1$ interpolation as a harmless simplification.

The focus of this work is twofold. First, we aim to provide a concise explanation of the most important analytical reasons that $H^1$ elements should be avoided in favor of $H(\Curl)$ edge elements in any $a$-$\phi$ formulation of electromagnetics. Second, we provide a concrete example in which $H^1$ elements fail to furnish a reliable approximation to the solution of a boundary value problem compared to $H(\Curl)$ elements, thereby making the issue concrete even to readers for whom the full analytical explanation may be difficult to parse. In the example, the quality of the solution is obviously very poor when fully $H^1$ elements are used, an issue that persists with refinement, while $H(\Curl)$ elements converge towards the correct solution. Through these explanations and the associated demonstrations, we illustrate that $H^1$ elements are structurally deficient in their ability to approximate solutions to problems in $a$-$\phi$ electromagnetics formulations compared to $H(\Curl)$ elements.

\section{Theoretical Background}

\subsection{Maxwell's equations}

The electromagnetic balance laws are widely summarized in many works \citep[see, e.g.,][among others]{kovetz:00}. In media, Maxwell's equations read as
\begin{align}
 \text{Gauss Law:}&&   \oint_{\partial\mathscr{B}}\bm{d}\cdot\bm{n}\,da&=\int_{\mathscr{B}}q\,dv \label{eq:me1}\\
  \text{Maxwell-Amp\`ere Law:}&&  \oint_{\partial\mathscr{A}}\bm{h}\cdot d\bm{x}&=\int_{\mathscr{A}}\bm{i}\cdot\bm{n}\, da+\frac{\partial}{\partial t}\int_{\mathscr{A}}\bm{d}\cdot\bm{n}\,da\label{eq:me2}\\
  \text{Faraday Law:}&&  \oint_{\partial\mathscr{A}}\bm{e}\cdot d\bm{x}&=-\frac{\partial}{\partial t}\int_{\mathscr{A}}\bm{b}\cdot \bm{n}\,da\label{eq:me3}\\
   \text{Monopole Law:}&& \oint_{\partial\mathscr{B}}\bm{b}\cdot\bm{n}\,da&=0\label{eq:me4},
\end{align}
where $\bm{d}$ is the electric displacement field\footnote{Also called the charge potential.} [C/m$^2$], $q$ is the macroscopic free charge density [C/m$^3$], $\bm{h}$ is the magnetic field\footnote{Also called the current potential.} [A/m], $\bm{i}$ is the free current density [A/m$^2$], $\bm{e}$ is the electric field [V/m], and $\bm{b}$ is the magnetic flux density\footnote{The term ``magnetic field'' is often used to refer to both $\bm{h}$ and $\bm{b}$, with meaning clear from context.} [T]. Here, $\mathscr{B}$ and $\mathscr{A}$ are arbitrary volumes and surfaces with unit normal vector fields $\bm{n}$, taken outward on $\partial\mathscr{B}$. It is worth noting that much of the difficulty associated with the interpolation of electromagnetic quantities is that as expressed in \eqref{eq:me1}-\eqref{eq:me4}, the fields are represented by their vector proxies, whereas 
in their ``true definition'' they are differential forms which respect a strict differential-geometric structure (see, e.g., \citet{Bossavit1998-GeomEM-1, Bossavit1998-GeomEM-2,Bossavit1998-GeomEM-3,Bossavit1998-GeomEM-4}).  The existence of the vector proxies is a consequence of the special structure of $\mathbb{R}^4$ (space-time),\footnote{The vector proxies are understood relative to an observer-dependent $3+1$ split, using the induced spatial metric and orientation.} and while the vector proxies make electromagnetics more accessible, they hide essential and important geometric structure of the governing equations and the fields.

In a stationary linear isotropic dielectric or magnetic material and a linear isotropic ohmic conductor, to close the system of equations, \eqref{eq:me1}-\eqref{eq:me4} are supplemented with the following constitutive equations:
\begin{align}
    \bm{d}=\epsilon_0(1+\chi)\bm{e}, && \bm{h}={\mu}_0^{-1}(1-\chi_m)\bm{b}, && \bm{i}=\kappa\bm{e}\ \label{eq:emfield}.
\end{align}
These are simplifications of more general models derivable from the electromagnetic duality relations and thermodynamic arguments \citep[see e.g.,][]{CARLSON:26}. Here $\epsilon_0$, ${\mu}_0$, $\kappa$, $\chi$, and $\chi_m$ are the (electric) permittivity of free space $\approx 8.854\times 10^{-12}$ [F/m], (magnetic) permeability of free space $\approx4\pi\times 10^{-7}$ [T$\cdot$m/A], electric conductivity [S/m], dielectric susceptibility, and magnetic susceptibility, respectively.  Also worth noting is that ${(\mu_0\epsilon_0)}^{-1} = c^2$, the speed of light squared.

\subsection{Time domain $a$-$\phi$ formulation}

On a contractible space (e.g., $\mathbb{R}^n$), the Poincar\'e lemma \citep[see, e.g.,][Chap.\ 17]{lee:12} implies that Faraday's and the Monopole laws \eqref{eq:me3} and \eqref{eq:me4} are automatically satisfied by introducing potentials $\bm{a}$ [T$\cdot$ m] and $\phi$ [V] in terms of which one declares that
\begin{align}
    \bm{b}=\nabla\times \bm{a}, \quad\text{and}\quad \bm{e}=-\bm{a}_{,t}-\nabla\phi.\label{eq:emagpots}
\end{align}
Therefore, Maxwell's equations can be reduced locally to
\begin{align}
\cdiv[\epsilon_0(1+\chi)(\bm{a}_{,t}+\nabla\phi)]+q&=0\\
    {\mu}_0^{-1}\nabla\times[(1-\chi_m)\nabla\times\bm{a}]+\kappa(\bm{a}_{,t}+\nabla\phi)+\epsilon_0[(1+\chi)(\bm{a}_{,t}+\nabla\phi)]_{,t}&=\bm{0},
\end{align}
modulo a to-be-discussed gauge constraint, say, the Coulomb or Lorentz gauge \citep[][\S6.3]{jackson:21}.
In many time domain applications, the frequency $\omega$ of electromagnetic loading is such that $\omega\epsilon_0\ll\kappa$, making the third term in the Maxwell-Amp\`{e}re law
negligible in comparison to the second term.
Further, materials can be assumed to be electro-neutral at the macroscopic scale so that $q=0$ in most cases. Thus, in practical applications it is often the case that the equations
\begin{align}
\cdiv[\epsilon_0(1+\chi)(\bm{a}_{,t}+\nabla\phi)]&=0,\label{eq:gaussref}\\
    {\mu}_0^{-1}\nabla\times[(1-\chi_m)\nabla\times\bm{a}]+\kappa(\bm{a}_{,t}+\nabla\phi)&=\bm{0}
\end{align}
determine the electromagnetic fields to good accuracy. It is worth remarking that, in materials with piecewise constant properties, \eqref{eq:gaussref} is locally equivalent to the charge conservation equation
\begin{align}
    \cdiv[\kappa(\bm{a}_{,t}+\nabla\phi)]&=0.\label{eq:cc}
\end{align}
The choice of which equation to solve, \eqref{eq:gaussref} or \eqref{eq:cc}, typically depends on whether the materials of interest are conductors or dielectrics; for example, in an eddy current application \eqref{eq:cc} would be preferred. 

In this work, only the $a$-$\phi$ formulation will be considered given its ability to be extended to a very general class of problems. Moreover, an $a$-$\phi$ formulation is often  necessary in multiphysics applications involving electrochemistry, where the electromagnetic potentials  become physically important in constitutive relations rather than just mathematically convenient \citep[see, e.g.,][]{CARLSON:26}.

\subsection{$H^1$ and $H(\Curl)$ elements: The discrete de Rham complex}

The potential formulation of Maxwell's equations is valid in large part thanks to the exactness of the de Rham complex over $\mathbb{R}^3$ or any contractible smooth submanifold \citep[see, e.g.,][Chap.\ 17]{lee:12}. Letting $\Pi_i$ denote the projections onto the finite element interpolation spaces used, one should ensure that the diagram
\begin{center}
    \begin{tikzcd}
        H^1\arrow[r, "\nabla"]\arrow[d,"\Pi_\phi"] & H(\Curl) \arrow[r, "\nabla\times"]\arrow[d,"\Pi_{\bm{a}}"] &  H(\Div) \arrow[r, "\cdiv"]\arrow[d,"\Pi_2"] & L^2 \arrow[d,"\Pi_3"]\\
        S_\phi \arrow[r, "\nabla_d"] & N_{\bm{a}} \arrow[r, "\nabla_d\times"] &  W_{\bm{b}} \arrow[r, "\cdiv_d"] & D_q
    \end{tikzcd}
\end{center}
commutes in order that approximated solutions satisfy the geometric structure of Maxwell's equations \citep[see, e.g.,][]{Hiptmair:02}. Here, $S_\phi$ denotes the interpolation space for $\phi$, $N_{\bm{a}}$ denotes the interpolation space for $\bm{a}$, $W_{\bm{b}}$ denotes an interpolation space for 2-form proxies like $\bm{b}$, $D_q$ denotes an interpolation space for 3-form proxies like $q$, and the $d$ subscript denotes the action of the discrete differential operator.\footnote{These discrete operators often act on reference cells and must be consistently pushed forward and pulled back. Hence, the isoparametric maps used are not always obvious on each of the spaces above \citep[see, e.g.,][]{Zaglmayr:06}.} For a given domain $\mathscr{D}$, these spaces are subspaces of $H^1=\{v\in L^2(\mathscr{D}) : \nabla v\in L^2(\mathscr{D})^3\}$, $H(\Curl)=\{\bm{v}\in L^2(\mathscr{D})^3 : \nabla\times\bm{v}\in L^2(\mathscr{D})^3\}$, $H(\Div)=\{\bm{v}\in L^2(\mathscr{D})^3 : \cdiv{\bm{v}}\in L^2(\mathscr{D})\}$, and $L^2=\{v\text{ measurable} : \int_{\mathscr{D}}|v|^2\,dx<\infty\}$, respectively, with derivatives understood weakly.

By commutation we mean the interpolation of the gradient of a scalar function should always be equivalent to the discrete gradient of the interpolation of that scalar function, and so on. Interpreted physically, the issue of interpolation is not merely whether the polynomial order is high enough, but whether the discrete spaces preserve the differential relationships among potentials, fields, fluxes, and charges.

A primary advantage of $H(\Curl)$ edge elements is their role in such a discrete complex, enabling good approximations of boundary value problem solutions to the order of approximation. On top of this, assignment of boundary conditions is far simpler and normal continuity of the time derivative of the vector potential need not be maintained; such discontinuities are often necessary for consistent solutions \citep[see, e.g.,][]{Jin:14}. On the other hand, using $H^1$ elements to approximate, for example, vector (1-form) objects like $\bm{a}$ is problematic in general given that these elements are not a part of any sequence satisfying this projection property. Therefore, satisfying Maxwell's equations in a weak, discrete setting with $H^1$ elements does not imply that the finite element solution is an approximation of the continuous, analytical problem in general, especially when the true solution is not $H^1$ continuous. 

\section{Numerical solution strategies}

\subsection{Weak form equations}

The equations to be solved in the present context over an approximation domain $\mathscr{D}$ can be stated in weak form as
\begin{align}
    \int_\mathscr{D}\kappa(\bm{a}_{,t}+\nabla\phi)\cdot\nabla{\delta \phi}\,dv+\int_{\partial\mathscr{D}}\delta \phi\,\bm{i}\cdot\bm{n}\,da=0,\label{eq:wfGcharge}
\end{align}
\begin{align}
    \int_\mathscr{D}(1-\chi_m)[\nabla\times\bm{a}]\cdot[\nabla\times\delta\bm{a}]&+\mu_0\kappa(\bm{a}_{,t}+\nabla\phi)\cdot\delta\bm{a}\,dv \nonumber \\ 
    &+\int_{\partial\mathscr{D}}(1-\chi_m)(\bm{n}\times\bm{b})\cdot\delta\bm{a}\,da=0,\label{eq:wfGvpot}
\end{align}
where $\delta\phi$ and $\delta\bm{a}$ are test functions. In general, the solution to this set of equations requires an additional specification of a gauge controlling the gradient null-space of $\bm{a}$ and the time dependent null-space of $\phi$. Heuristically, these are the ``rigid body'' modes; if $\hat{\bm{a}}$ and $\hat{\phi}$ are solutions, so are $\hat{\bm{a}}+\nabla f$ and $\hat{\phi}-f_{,t}$ for any scalar field $f$. A common and simple choice is the Coulomb gauge condition, $\cdiv\bm{a}=0$ \citep[][\S6.3]{jackson:21}.

\subsubsection{Edge element formulations}\label{s:edge_form}

In this work, mixed edge scalar $H^1/H(\Curl)$ elements, $H^1$ for $\phi$ and $H(\Curl)$ for $\bm{a}$, are used to interpolate solutions to boundary value problems. The core advantage of edge elements is that the associated interpolation space $N_{\bm{a}}$ is a part of a discrete exact complex. 
For the transient problems presently considered, the time-discretized form of Eqs.~\eqref{eq:wfGcharge}--\eqref{eq:wfGvpot} is coercive on the quotient space of ordered pairs $(\bm{a},\phi)$ modulo gauges ($(\nabla f,-f_{,t})$ for any admissible scalar field $f$), in the quotient norm induced by the mixed $H(\Curl)\times H^1$ norm.\footnote{Equivalently, in the continuum, space-time 1-forms modulo exact space-time 1-forms.} Hence, for compatible data, correct solutions can be obtained via any iterative method that converges to a representative of the solution class without the need to explicitly impose a gauge, with the unknown gauge being of no consequence to the physical quantities of interest.
The commuting discretization preserves the required compatibility structure at the discrete level, allowing the gauge to be chosen at each timestep while leaving the computed electromagnetic fields invariant. This is analogous to the invariance of material stress and strain under superposed rigid body motion in elasticity. 

\begin{figure}\centering
\includegraphics[height=0.25\textheight]{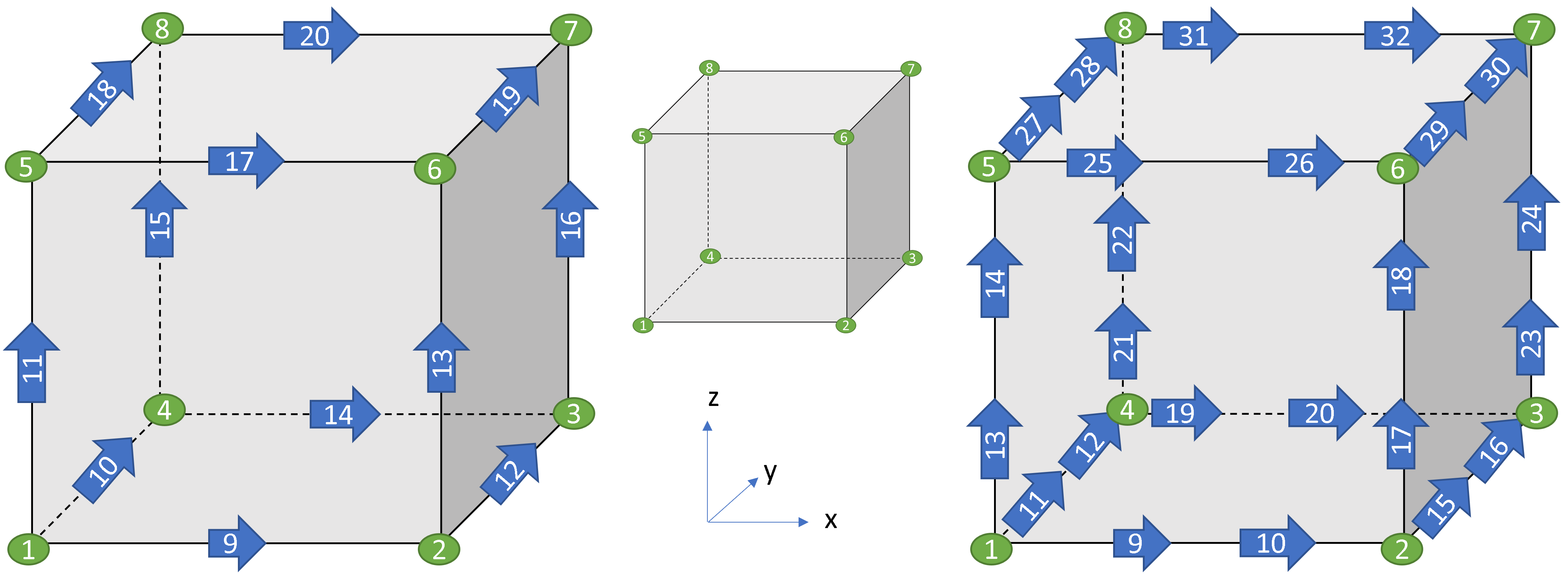}
\caption{Depictions of the 20 (left) and 32 (right) node elements utilized in this work. The typical 8 node hexahedral reference element is also shown (middle).} \label{f:1}
\end{figure}

The mixed edge scalar elements implemented in this study are 20 and 32 node hexahedral elements, as depicted in Fig.~\ref{f:1}. The 20 node elements consist of eight corner nodes associated with first order, trilinear degrees of freedom used to interpolate the scalar potential $\phi$ and a single node on each edge associated with zeroth-order complete N\'ed\'elec edge degrees of freedom used to interpolate the magnetic vector potential $\bm{a}$. These two spaces form a part of a discrete exact complex, meaning that $\nabla S_\phi=\ker(\nabla\times|_{N_{\bm{a}}})$, mirroring the exact continuous complex.

In the time domain, exactness of the discrete complex is not strictly necessary to furnish good solutions to boundary value problems with mixed scalar/edge interpolations so long as the exact complex associated with the lowest order space is a subcomplex of the exact complex associated with the highest order space. The commutative nature of the lowest order diagram ensures that the higher order quantity is projected consistently in the lowest order equation. Thus, any error associated with this particular form of inexactness is identical in nature to the error associated with the general interpolation of, for example, a nonlinear solution with linear elements, and is not a more fundamental structural issue leading to genuinely spurious solutions at lowest order. To demonstrate how these elements can be effective, the 32 node elements consist of the same scalar corner degrees of freedom as in the 20 node case, only now with 24 first-order complete serendipity edge degrees of freedom for the magnetic vector potential with two degrees of freedom associated with each edge \citep{arndold:14}. These serendipity $H(\Curl)$ functions are a part of a discrete exact complex containing the complex associated with the first-order $H^1$ scalar degrees of freedom as a subcomplex. Again, although the discrete complex formed by $S_\phi$ and $N_{\bm{a}}$ is now not exact, the continuous geometric structure is still maintained to lowest order while simultaneously allowing for a higher order interpolation of a degree of freedom of choice: in this case, the vector potential $\bm{a}$. Thus, one can achieve a higher order of accuracy with a reduced number of degrees of freedom; in this case, the first-order $H(\Curl)$ elements have the same number of degrees of freedom as the fully $H^1$ elements.\footnote{On assembly, a fully $H^1$ mesh will tend to have fewer degrees of freedom given that a typical node in a regular mesh is shared by eight elements, whereas an edge is only shared by four elements.}

\subsubsection{Fully $H^1$ element formulations}

As noted, fully $H^1$ elements where each component of $\bm{a}$ is interpolated as an $H^1$ continuous scalar are not a part of a discrete exact complex that forms a commutative diagram with the continuous de Rham sequence. Therefore, modes that discretely satisfy \eqref{eq:wfGcharge} and \eqref{eq:wfGvpot} may not correspond to actual, continuous fields associated with the boundary and initial data. In many cases, this failure is especially evident for the gradient part of $\bm{a}$ in the Helmholtz decomposition. To stabilize/control these spurious modes, a penalty term is often added to \eqref{eq:wfGvpot} in $H^1$ formulations to enforce the Coulomb gauge, leading to a modified weak form:
\begin{equation}\label{eq:Pform}\begin{aligned}
    \int_\mathscr{D}(1-\chi_m)\left([\nabla\times\bm{a}]\cdot[\nabla\times\delta\bm{a}]+\lambda\cdiv(\bm{a})\cdiv(\delta\bm{a})\right)+\mu_0\kappa(\bm{a}_{,t}+\nabla\phi)\cdot\delta\bm{a}\,dv
    \\+\int_{\partial\mathscr{D}}(1-\chi_m)(\bm{n}\times\bm{b})\cdot\delta\bm{a}\,da=0,
\end{aligned}\end{equation}
where $\lambda$ is a penalty parameter often set to unity (see, e.g., \cite{coulomb:81}, \citet[][\S5.7.4]{Jin:14}, or \cite{anand:26}).\footnote{Eqn.\ \eqref{eq:Pform} is typically scaled by $\mu_0^{-1}$, and thus the penalty is usually written as either $\mu_0^{-1}$ or $(1-\chi_m)\mu_0^{-1}$. Although different at first glance, the penalty parameter above agrees with the typical value used by, e.g., \citet{Stewart:26} after noting the scaling.} This gauge is generally necessary since the structure of the function spaces does not allow for consistent gauge evolution in time without explicit control. 

It should be noted that, even if strictly enforced, this constraint is not enough to ensure that a discretely computed curl corresponds to the compatible projection of the continuous curl associated with the boundary data. To fully recover the correct geometric structure, one should also solve $\cdiv{\bm{b}}=0$ with a separate interpolation for $\bm{b}$ under the constraint $\bm{b}=\nabla\times\bm{a}$. Although seemingly redundant, strictly speaking such compatibility is necessary given that, for arbitrary $\bm{a}'\in H(\Curl)$, it may be the case that $\Pi_2(\nabla\times\bm{a}')\neq\nabla_d\times(\Pi_{\bm{a}}\bm{a}')$. In other words, even a strict divergence constraint on $\bm{a}$ cannot guarantee the absence of spurious solutions given that a discretely computed curl need not correspond to the continuous curl of interest.

The $H^1$ elements implemented in this work are 32 degree of freedom eight node trilinear hexahedral elements consisting of four degrees of freedom at each node: three to interpolate each component of the magnetic vector potential $\bm{a}$ and one to interpolate the electric scalar potential. This element is graphically depicted in Fig.~\ref{f:1}.

\subsection{Solver technology}

The linearized system of equations for problems posed by the theory is ill conditioned, even with spatially uniform
conductivity values. To still provide reasonable solutions at an acceptable tolerance,\footnote{The default FEAP tolerance is used in this work: $10^{-8}$ relative residual norm and $10^{-16}$ relative energy norm; see \cite{taylor_feap:2014} for more information.} the field-split preconditioning methodology of \cite{CARLSON:26} is utilized in PETSc with a flexible, iterative GMRES solver \citep{petsc-efficient,petsc-user-ref}. This solver, along with the elements discussed, is implemented in FEAP \citep{taylor_feap:2014}. Time stepping is carried out using FEAP's BDF2 time integrator, and all numerical integration is computed via a 3$\times$3$\times$3 Gaussian quadrature rule. 

\section{Adversarial benchmark problem}

To illustrate the consequences of the  limitations of the $H^1$ approximation discussed above, a bi-material square prism subjected to a linearly increasing magnetic field is considered. The geometry and loading admit an independently constructed semi-analytic reference solution, allowing the current density computed by the finite element formulations to be assessed directly. After a short startup transient, the induced current becomes independent of time, providing a simple setting in which to examine approximation error in space and time. 

\subsection{Geometry, material properties, and loading}

The specimen occupies a square prism of dimensions $1.5\,\mathrm{mm}\times1.5\,\mathrm{mm}\times0.5\,\mathrm{mm}$, centered at the origin. An aluminum alloy conducting annulus surrounds a titanium alloy core with a $0.5\,\mathrm{mm}\times0.5\,\mathrm{mm}$ cross section extending through the full thickness of the specimen. Thus, the annular band has a thickness of $0.5\,\mathrm{mm}$ on each side. The material properties are idealized as indicated in Table~\ref{tab:square-materials}; in particular, the conductivity contrast is prescribed to be exactly ten to one.\footnote{The conductivity values given fall within the realizable ranges for both aluminum and titanium alloys.} Both materials are treated as non-magnetic, and perfect electrical contact is assumed at their interface.

\begin{table}[htbp]
  \centering
  \caption{Simulation model parameters}
  \label{tab:square-materials}
  \begin{tabular}{lcc}
    \toprule
    Parameter & Aluminum alloy annulus & Titanium alloy core \\
    \midrule
    Conductivity $\kappa$ [S/m]
      & $\kappa_a=2.538\times10^7$ & $\kappa_a/10=2.538\times10^6$ \\
    Magnetic susceptibility $\chi_m$ [$-$] & $0$ & $0$ \\
    \bottomrule
  \end{tabular}
\end{table}

\begin{figure}[htbp]
  \centering
  \IfFileExists{square_benchmark_geometry.pdf}{%
    \includegraphics[width=0.9\textwidth]{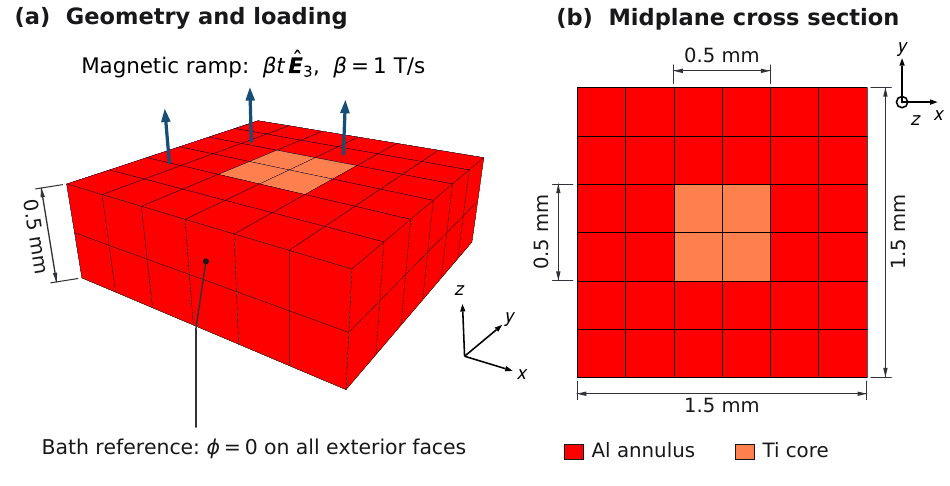}%
  }{%
    \fbox{\parbox[c][35mm][c]{0.72\textwidth}{\centering
      Insert the new square-prism geometry and mesh.\\[0.5em]
      Show the core, annulus, axial magnetic ramp, and $\phi=0$.}}%
  }
  \caption{Geometry and loading of the bi-material conductor. The titanium core extends through the full thickness of the aluminum annulus. On any given exterior face with normal $\bm{n}$ and at any time $t$, the applied magnetic boundary condition is $\bm{n}\times\bm{b}=\beta t(\bm{n}\times\hat{\bm{E}}_3)$.}
  \label{f:square-geometry}
\end{figure}

Denote the domain by $\mathscr{D}$ and its exterior boundary by $\partial\mathscr{D}$. The imposed magnetic boundary condition and scalar potential are\footnote{Note that $\hat{\bm{E}}_i$ denotes the $i$th coordinate direction.}
\begin{equation}
  \bm{n}\times\bm{b}
  =\bm{n}\times(\beta t\,\hat{\bm{E}}_3),
  \qquad \phi=0
  \quad\text{on }\partial\mathscr{D},
  \qquad \beta=1\,\mathrm{T/s}.
  \label{eq:square-loading}
\end{equation}
The ramp begins at $t=0$ from an initially field-free, current-free state, with $\bm{a}_0=\bm{0}$ and $\phi_0=0$. In the magneto-quasi-static limit of the Maxwell--Amp\`ere law, Eq.~\eqref{eq:me2}, the prescribed tangential magnetic field has zero circulation around any closed loop on an exterior face and therefore implies $\bm{i}\cdot\bm{n}=0$ on each such face. The geometry and boundary conditions are depicted in Fig.~\ref{f:square-geometry}.

The fixed scalar potential is interpreted as a prescribed electrochemical reference state. A possible physical idealization is a specimen equilibrated in a reversible redox bath, with the whole assembly placed in a nearly uniform, axially ramped solenoid $\bm{b}$-field. Given the small specimen size and electric field magnitudes, it is reasonable to assume that chemical equilibrium is approximately maintained at the specimen–bath interface and to neglect the microtesla-scale self-field when identifying the prescribed boundary field with the applied field.

\begin{remark}
    Finite element solutions are computed on three-dimensional meshes despite the two-dimensional nature of the eddy-current solution.
\end{remark}

\subsection{Semi-analytic reference solution}

By symmetry of the geometry and loading, the steady-ramp magnetic field
can be written as $\bm b=(\beta t-g)\hat{\bm{E}}_3$, where $g=g(x,y)$ vanishes on the perimeter of the cross section. The constitutive relations in
Eq.~\eqref{eq:emfield}, together with the magneto-quasi-static
Maxwell--Amp\`ere and Faraday laws, Eqs.~\eqref{eq:me2} and~\eqref{eq:me3},
reduce the problem to
\begin{equation}
  -\cdiv(\kappa^{-1}\nabla g)=\mu_0\beta .
  \label{eq:square-poisson}
\end{equation}
Because both materials are non-magnetic and perfect electrical contact is assumed, $g$ and $\kappa^{-1}\nabla g\cdot\bm{n}$ are continuous across material interfaces.\footnote{Continuity of $\kappa^{-1}\nabla g\cdot\bm{n}$ across material
interfaces follows from tangential continuity of $\bm{e}$. Indeed, by
Eq.~\eqref{eq:square-reference-fields} and Eq.~\eqref{eq:emfield},
$\kappa^{-1}\nabla g\cdot\bm{n}
=\frac{\mu_0}{\kappa}(\bm{n}\times\bm{i})\cdot\hat{\bm{E}}_3
=\mu_0\bm{e}\cdot(\hat{\bm{E}}_3\times\bm{n})$,
where $\hat{\bm{E}}_3\times\bm{n}$ is tangent to the interface.} The resulting eddy current is
\begin{equation}
  \bm i_{\mathrm{ref}}=\frac{1}{\mu_0}\left(\frac{\partial g}{\partial x}\hat{\bm{E}}_2-\frac{\partial g}{\partial y} \hat{\bm{E}}_1\right).
  \label{eq:square-reference-fields}
\end{equation}

Details regarding the semi-analytic solution for \(g\) are given in Appendix~\ref{app:square-reference}.

\section{Results}
\label{sec:square-results}

In all finite element calculations, the volume elements are cubes and the time step is fixed at $\Delta t=5\,\mu\mathrm{s}$. Computations are carried out on a single processor and, unless otherwise specified, the penalty parameter is set to $\lambda=1$ in $H^1$ element formulations.

\subsection{Solution comparisons}

\begin{figure}[htbp]
  \centering
  \IfFileExists{square_benchmark_comparison.pdf}{%
    \includegraphics[width=\linewidth]{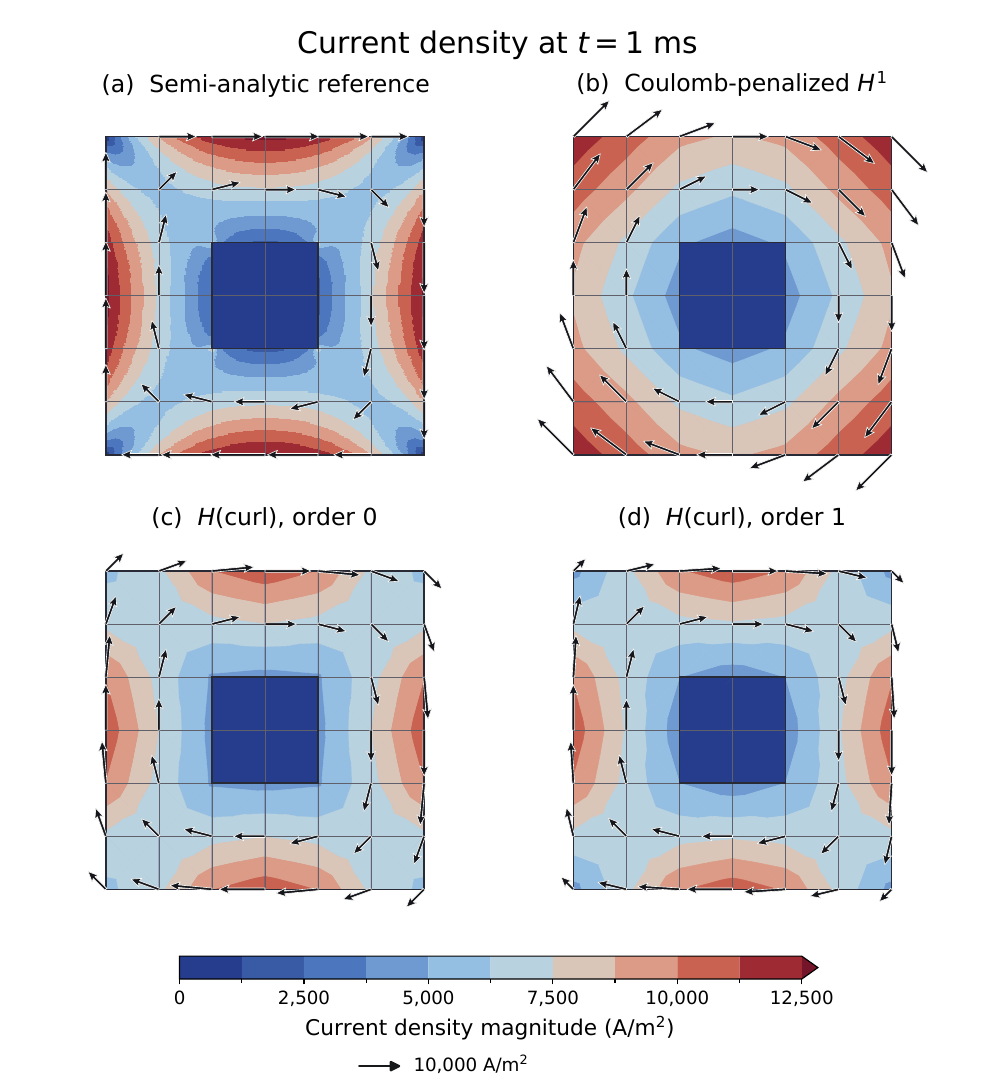}%
  }{%
    \fbox{\parbox[c][40mm][c]{0.90\linewidth}{\centering
      Coarse-mesh current fields and semi-analytic reference}}%
  }
  \caption{Eddy-current distributions on the coarsest mesh under a linear magnetic-field ramp. The finite element fields and the semi-analytic reference are shown with identical limits for current magnitude and a common scale for the current vectors. The finite element mesh is superimposed on the semi-analytic reference field solely to aid comparison.}
  \label{f:square-comparison}
\end{figure}

Fig.~\ref{f:square-comparison} compares the computed current fields on the coarsest mesh at $t=1\,\mathrm{ms}$ with the semi-analytic steady-ramp solution of Eq.~\eqref{eq:square-reference-fields}. In the reference solution, the current circulates through the aluminum annulus, follows the exterior faces, and vanishes at the four outer corners of the cross section. However, the Coulomb-penalized $H^1$ solution places its largest currents at precisely these corners. Additionally, along the adjacent faces, the current vectors exhibit a pronounced component normal to the exterior boundary, in violation of the no-through-current condition $\bm{i}\cdot\bm{n}=0$ implied by the magnetic field boundary condition. The $H(\Curl)$ solutions more closely reproduce both the distribution of current in the annulus and its direction near the boundary. Moreover, for both $H(\Curl)$ formulations, the settled eddy-current profile approaches the reference profile as the mesh is refined (Fig.~\ref{f:square-error-history}). In contrast, on all three meshes the settled Coulomb-penalized $H^1$ solution remains similar to the coarse-mesh solution shown in Fig.~\ref{f:square-comparison}(b).

\subsection{Effect of mesh refinement}

\begin{figure}[htbp]
  \centering
  \includegraphics[width=\linewidth]{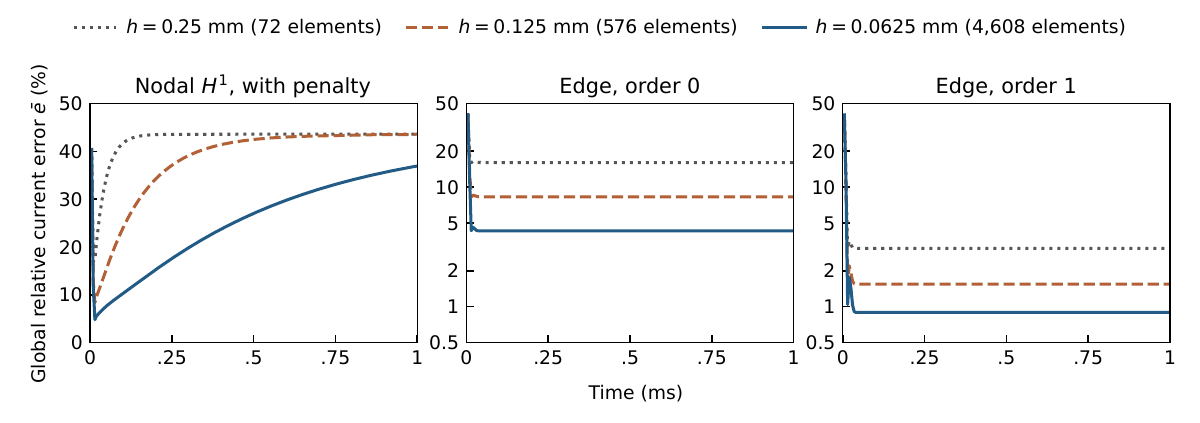}
  \caption{Relative current error for the Coulomb-penalized $H^1$, the order-zero $H(\Curl)$, and the order-one $H(\Curl)$ element formulations. Here, $h$ denotes the element edge length, and the middle and right panels use semi-log axes.}
  \label{f:square-error-history}
\end{figure}

Fig.~\ref{f:square-error-history} compares three meshes containing 72, 576, and 4608 volume elements, obtained by successively halving the element edge length. The relative current error,
\begin{equation}
  \bar{e}(t)=
  \frac{\|\bm{i}_h(t)-\bm{i}_{\mathrm{ref}}\|_{\kappa^{-1}}}
       {\|\bm{i}_{\mathrm{ref}}\|_{\kappa^{-1}}},
  \label{eq:square-error-measures}
\end{equation}
is evaluated from the element fields at volume quadrature points using the resistivity-weighted norm 
\begin{equation}
  \|\bm i\|_{\kappa^{-1}}^2=\int_{\mathscr D}\kappa^{-1}|\bm i|^2\,dv .
  \label{eq:square-current-norm}
\end{equation}
Because the reference is the steady-ramp solution, $\bar{e}$ initially includes the physical startup transient. However, according to Eq.~\eqref{eq:square-transient-time} the continuum current discrepancy has fallen below $0.1\%$ of its initial value by $25.2\,\mu\mathrm{s}$.

For both $H(\Curl)$ formulations, the error settles within the first ten increments ($50\,\mu\mathrm{s}$) to a level that decreases with each refinement. At $t=0.15\,\mathrm{ms}$, the order-zero errors are $16.07\%$, $8.31\%$, and $4.31\%$, from coarsest to finest. The corresponding order-one errors are $3.07\%$, $1.55\%$, and $0.895\%$. Both $H(\Curl)$ formulations correctly predict that the eddy-current field becomes time-independent after the startup transient, even as the prescribed magnetic field continues to increase linearly, as in Eq.~\eqref{eq:square-reference-fields}. 

By contrast, the penalized $H^1$ current continues to evolve well beyond the physical startup transient and ultimately settles to an incorrect current distribution (Fig.~\ref{f:square-comparison}(b)). Mesh refinement initially reduces the current error but increases the duration of this erroneous transient. As shown in Fig.~\ref{f:square-error-history}, the error rises back above $30\%$ at $0.040$, $0.155$, and $0.610\,\mathrm{ms}$, respectively, from coarsest to finest. Thus, halving the element size delays this rise by approximately a factor of four at each refinement. At $1\,\mathrm{ms}$, the errors on the two coarser meshes are nearly identical ($43.60\%$ and $43.65\%$), whereas the finest-mesh error is lower ($36.95\%$) and still rising. Additional time increments confirm that the error on the finest mesh eventually reaches approximately the same plateau as on the coarser meshes. Refinement therefore delays the rise in error but does not correct the settled current distribution.

\section{Discussion}
\label{sec:square-discussion}

Additional insight is obtained by comparing the present results with those for the TEAM 4 problem \citep{Kameari:88} and by considering the penalty-free order-one $H^1$ formulation. In the TEAM 4 problem, all element formulations considered in this work produce reasonable eddy-current solutions, including the penalty-free ($\lambda=0$) $H^1$ formulation. By contrast, neither the penalized nor the penalty-free $H^1$ formulation reproduces the settled current distribution in the present benchmark. This suggests that the TEAM 4 problem does not strongly excite the spurious discrete gradient components identified in the preceding analysis, whereas the present benchmark is sensitive to them. These components belong to $\bigl[\Pi_{\bm{a}}(H(\Curl))/\Pi_{\bm{a}}(\ker(\nabla\times))\bigr]\cap\nabla_d S_\phi$,\footnote{The intersection is understood in the quotient space; its nonzero classes are represented by discrete gradients that are not projections of continuous curl-free fields.} which is nontrivial for the $H^1$ approximation considered here. The TEAM 4 results raise the possibility that many standard eddy-current benchmarks may not strongly excite this space. Agreement on such benchmarks alone could therefore give the impression that nodal elements generally entail only a minor loss of accuracy. In the present benchmark, however, the discrepancy persists under sustained loading and is apparent not only in the current magnitude but also in its direction at the exterior boundary (Fig.~\ref{f:square-comparison}). The pronounced normal component of the computed current is consistent with the loss of the discrete de Rham structure discussed previously: the continuum condition $\bm{i}\cdot\bm{n}=0$ follows from the prescribed tangential magnetic field in Eq.~\eqref{eq:square-loading}, without being imposed separately in the $\phi$ equation.

\begin{figure}[htbp]
  \centering
  \includegraphics[width=\linewidth]{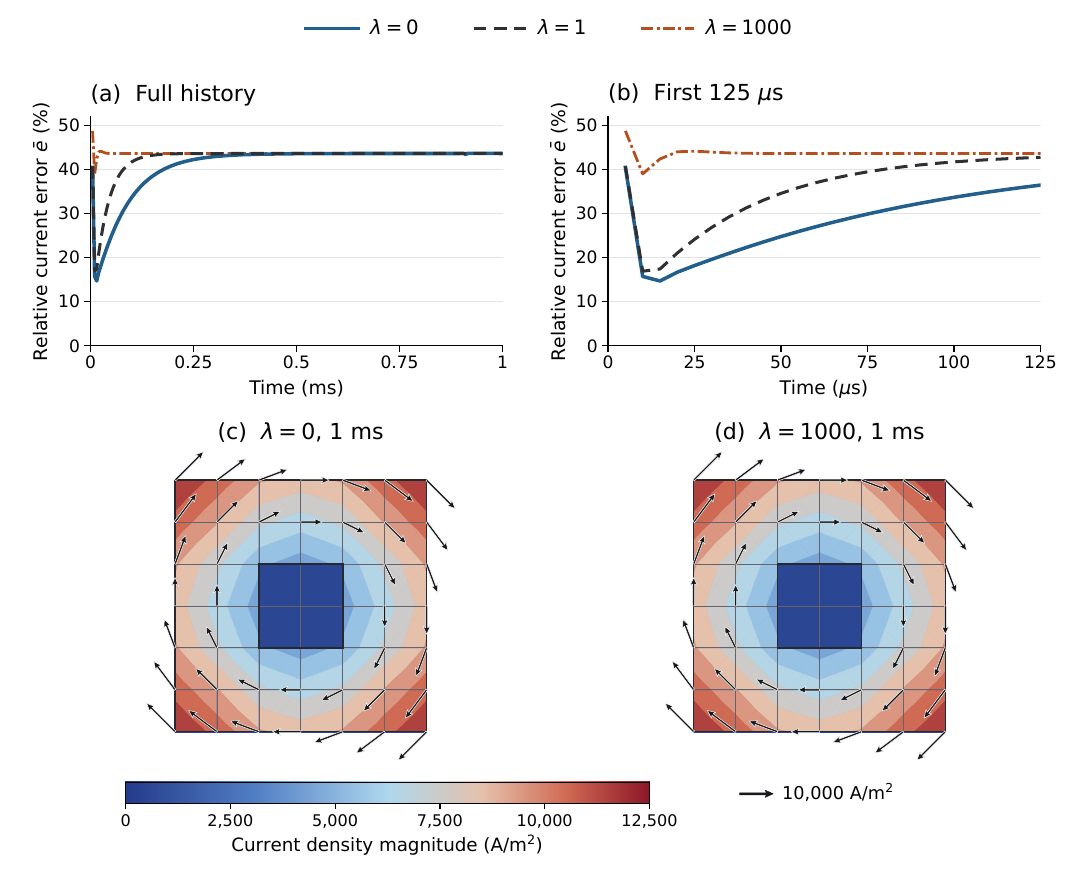}
  \caption{Effect of the Coulomb penalty on the order-one nodal $H^1$
    formulation on the 72-element mesh.
    (a) Relative current error of Eq.~\eqref{eq:square-error-measures}
    at all 200 recorded increments.
    (b) The first $125\,\mu\mathrm{s}$ of the same histories.
    (c,d) Top-face current-density magnitude and in-plane current vectors at
    $1\,\mathrm{ms}$ for the penalty-free and increased-penalty cases,
    respectively, with common magnitude and arrow scales.}
  \label{f:square-penalty-sensitivity}
\end{figure}

Moreover, it is evident from Fig.~\ref{f:square-penalty-sensitivity} that one cannot simply increase or remove the penalty parameter to improve the accuracy of the settled $H^1$ solution. On the 72-element mesh, removing the Coulomb gauge penalty (i.e., setting $\lambda=0$) delays the rise back above $40\%$ error from $0.080$ to $0.180\,\mathrm{ms}$. With $\lambda=1000$, which is 1000 times the conventional value, the error instead never drops below $39\%$.  
After the decay of the startup transient, the error rises almost immediately above $44\%$ and then approaches the plateau of $43.6\%$ from above, whereas the other two $H^1$ solutions approach it from below. The current plots are likewise strikingly similar, with each $H^1$ solution exhibiting a pronounced current component normal to the exterior faces. This example thus demonstrates that Coulomb gauge enforcement does not correct the underlying incompatibility of the approximation spaces, and that the failure to satisfy the implied condition $\bm{i}\cdot\bm{n}=0$ persists with or without the penalty term.

Finally, it is worth noting that the order-one $H(\Curl)$ formulation on the coarsest mesh gives a smaller relative current error than the order-zero formulation on any of the meshes considered (Fig.~\ref{f:square-error-history}). At $t=0.15\,\mathrm{ms}$, the error on the 72-element mesh is $3.07\%$ for the order-one formulation, compared with $16.07\%$ for the order-zero formulation. Even on the 4608-element mesh, the order-zero error is $4.31\%$. For the same 30 time increments, the coarse order-one solution required only $4.14$ CPU seconds, compared with $1237.40$ CPU seconds for the finest order-zero solution. Thus, by the chosen error measure, the more accurate solution was obtained in approximately one three-hundredth of the computational time. This result is of particular interest because only the vector potential was enriched: the scalar potential retained its trilinear interpolation. Although the resulting mixed $H^1/H(\Curl)$ pair is not exact, it preserves the lower-order discrete subcomplex. The present example thus emphasizes that, in time-domain problems of this type, preserving this subcomplex may be more important for current accuracy than requiring the mixed pair itself to be exact, potentially allowing for more efficient and accurate computation of solutions when one field requires a richer approximation than the other.

\section{Conclusion}

In this work, the time-domain $a$-$\phi$ electromagnetics equations and associated finite element implementations have been concisely summarized, and the issues plaguing $H^1$ continuous vector interpolations of the magnetic vector potential have been stated in their simplest form. This explanation is useful given that many similar explanations can be esoteric and difficult to quickly parse. Of primary importance is the clear demonstration of the failure of $H^1$ continuous vector-potential interpolations with penalty regularization of the (Coulomb) gauge in a simple benchmark problem. Although the $H^1$ solution initially appears to approach the semi-analytic reference, it ultimately settles to an incorrect current distribution that violates the no-through-current condition. Mesh refinement delays the approach to this state but does not correct it, nor does increasing or removing the penalty. In contrast, both $H(\Curl)$ formulations predict the expected time-independent eddy-current response within the first ten increments, and their settled current errors decrease with refinement. Moreover, the order-one formulation on the coarsest mesh gives a smaller relative current error than the order-zero formulation on the finest mesh in approximately one three-hundredth of the CPU time. Although this example demonstrates that exactness of the full mixed pair is not strictly necessary for an accurate current, both $H(\Curl)$ formulations preserve an exact lower-order discrete subcomplex, whereas the nodal $H^1$ formulation lacks this structure and fails to recover the correct settled current even under mesh refinement.

\bibliographystyle{elsarticle-harv}
\bibliography{used_references.bib}

\newpage
\backmatter
\begin{appendices}

\section{Construction of the reference solution}
\label{app:square-reference}

Coordinates in the series and edge integrals below are given in millimeters, with all other units in SI.

\subsection{Analytic solution form}

\begin{figure}[htbp]
  \centering
  \includegraphics[width=0.50\linewidth]{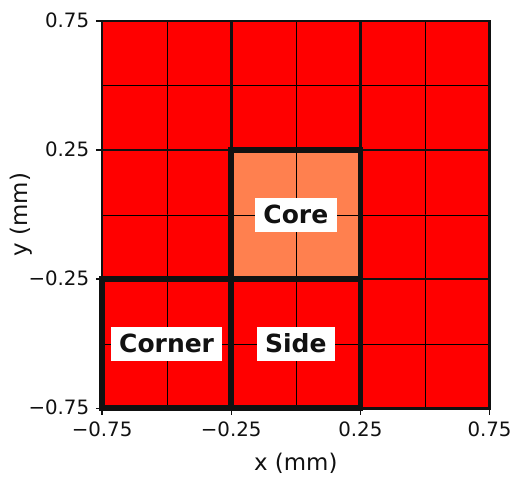}
  \caption{Partition of the cross section superimposed on the finite element mesh.}
  \label{f:square-reference-tiles}
\end{figure}

To solve Eq.~\eqref{eq:square-poisson}, extend the four sides of the Ti
core to the perimeter to create nine square subdomains of equal size.
By symmetry, only the three squares outlined in
Fig.~\ref{f:square-reference-tiles} need to be explicitly considered.
Denote the lower-left corner of the square being evaluated by $(x_c,y_c)$.
Its coordinates are $(-0.75,-0.75)$, $(-0.25,-0.75)$,
and $(-0.25,-0.25)$ for the corner, side, and core, respectively, with
$x=y=0$ at the center of the full cross section.

On each square, write $g=p_\kappa+h$, where $p_\kappa$ solves
$-\cdiv{\nabla p_\kappa}=\mu_0\beta\kappa$ with $p_\kappa=0$ on all four
sides. The harmonic correction satisfies $\cdiv{\nabla h}=0$ with $h=g$
on the sides of the square. The particular solution is the same on all
three squares up to a translation and the value of $\kappa$, and is given below following \citet[Secs.~2.5.1 and~8.6]{Haberman:04}:
\begin{equation}
 p_\kappa(x,y)=10^{-6}\mu_0\beta\kappa
 \sum_{\substack{n\geq1\\ n\ {\rm odd}}}
 \frac{\sin\!\big(2n\pi[x-x_c]\big)}{(n\pi)^3}
 \left[1-\frac{\cosh\!\big(2n\pi[y-y_c-0.25]\big)}
                  {\cosh(n\pi/2)}\right].
 \label{eq:square-particular}
\end{equation}

In addition, $h$ can be written as the bilinear interpolant $B(x,y)$ of the corner values of
$g$ plus a weighted sum of edge modes. By symmetry, $B=g=g_c$ at each
corner of the core, and from the applied boundary conditions $B=g=0$ at each corner on the outer perimeter.
The four edge modes for each $n$ are rotations of the top-edge mode
up to orientation:
\begin{equation}
\begin{aligned}
 H_n^{\rm top}(x,y)&=
 \sin\!\big(2n\pi[x-x_c]\big)\,\frac{\sinh\!\big(2n\pi[y-y_c]\big)}{\sinh(n\pi)},\\
 H_n^{\rm bottom}(x,y)&=
 \sin\!\big(2n\pi[x-x_c]\big)\,\frac{\sinh\!\big(2n\pi[y_c+0.5-y]\big)}{\sinh(n\pi)},\\
 H_n^{\rm right}(x,y)&=
 \sin\!\big(2n\pi[y-y_c]\big)\,\frac{\sinh\!\big(2n\pi[x-x_c]\big)}{\sinh(n\pi)},\\
 H_n^{\rm left}(x,y)&=
 \sin\!\big(2n\pi[y-y_c]\big)\,\frac{\sinh\!\big(2n\pi[x_c+0.5-x]\big)}{\sinh(n\pi)},
\end{aligned}
 \label{eq:square-edge-mode}
\end{equation}
where the edge coefficients are zero on the outer perimeter. Let $a_n$ denote the coefficients on the Al--Al cuts and $b_n$ those
on the Al--Ti interfaces. On each Al--Ti interface, $g-g_c$ is symmetric
about the midpoint, whereas the even-indexed sine modes are antisymmetric.
Thus, $b_n=0$ for even $n$, and the analytic solution on each square subset is
\begin{align}
 g_{\rm corner}(x,y)
 &=p_{\kappa_a}+4g_c(x-x_c)(y-y_c)
   +\lim_{N\to\infty}\sum_{n=1}^N a_n(H_n^{\rm top}+H_n^{\rm right}),
   \label{eq:square-tile-solution}\\
 g_{\rm side}(x,y)
 &=p_{\kappa_a}+2g_c(y-y_c)
   +\lim_{N\to\infty}\sum_{n=1}^N\big[a_n(H_n^{\rm left}+H_n^{\rm right})+b_nH_n^{\rm top}\big],
   \label{eq:square-side-solution}\\
 g_{\rm core}(x,y)
 &=p_{\kappa_a/10}+g_c
   +\lim_{N\to\infty}\sum_{n=1}^N b_n(H_n^{\rm bottom}+H_n^{\rm top}+H_n^{\rm left}+H_n^{\rm right}).
   \label{eq:square-core-solution}
\end{align}
Here, the remaining unknowns $g_c$, $a_n$, and the odd-indexed $b_n$ can be determined from the required tangential continuity of $\bm{e}$, which manifests as normal continuity of $\kappa^{-1}\nabla g$.

\subsection{Determination of the coefficients}

The tile solutions in Eqs.~\eqref{eq:square-tile-solution}--\eqref{eq:square-core-solution}
agree in value on each shared edge by construction, so it remains to match
$\kappa^{-1}\nabla g\cdot\bm{n}$. By symmetry, it suffices to consider
one of the eight Al--Al cuts and one of the four Al--Ti interfaces.
For the squares outlined in Fig.~\ref{f:square-reference-tiles}, these are chosen to
lie on $x=-0.25$ and $y=-0.25$, respectively. Let $w$ denote a scalar $H^1$
function defined on the full cross section that vanishes on the outer perimeter
and has the same symmetries as $g$. For every such $w$, the interface condition can be written in weak form as
\begin{equation}
\begin{aligned}
 &8\int_{-0.75}^{-0.25}
 \left.\left(\frac{\partial g_{\rm corner}}{\partial x}
            -\frac{\partial g_{\rm side}}{\partial x}\right)\right|_{x=-0.25}
             w(-0.25,y)\,dy\\
 &\quad+4\int_{-0.25}^{0.25}
 \left.\left(\frac{\partial g_{\rm side}}{\partial y}
          -10\frac{\partial g_{\rm core}}{\partial y}\right)\right|_{y=-0.25}
             w(x,-0.25)\,dx=0.
\end{aligned}
 \label{eq:square-flux-match}
\end{equation}
Here, the factor of ten accounts
for the lower conductivity of the core.

In order to compute the reference solution, the edge series is truncated at $n=N$.
The test space is chosen in a Bubnov--Galerkin sense: it is spanned by the functions multiplying $g_c$, $a_n$, and
$b_n$ in the edge values of $g$. Written as pairs $(w(-0.25,y),w(x,-0.25))$, these are
\[
 \big(2(y+0.75),1\big),\qquad
 \big(\sin(2n\pi[y+0.75]),0\big),\qquad
 \big(0,\sin(2n\pi[x+0.25])\big),
\]
respectively, with $1\le n\le N$ for $a_n$ and only odd $n\le N$ for $b_n$.
Evaluating the integrals in Eq.~\eqref{eq:square-flux-match} analytically for each test function gives one linear equation per coefficient, with the known terms from $p_\kappa$ supplying the right-hand side. The resulting symmetry-reduced system is then solved by LU factorization using \texttt{numpy.linalg.solve}.

To compute the reference solution used in Secs.~\ref{sec:square-results}--\ref{sec:square-discussion}, the coefficient system is solved with $N=512$, while $p_\kappa$ is evaluated separately according to Eq.~\eqref{eq:square-particular} by summing the purely sinusoidal part exactly and retaining the odd modes through $n=2047$ in the hyperbolic part. With this cutoff to $p_\kappa$ held fixed, doubling $N$ from 256 to 512 changes the current by only $0.0685\%$ relative to the 512-mode solution in the norm of Eq.~\eqref{eq:square-current-norm}. The reference current is therefore taken to be converged for the comparisons presented in this work.

\subsection{Comparison with the ramp started from zero}

The simulations start from zero current, whereas the steady-ramp current
in Eq.~\eqref{eq:square-reference-fields} is independent of time. Letting $r$ denote the difference between the sudden-ramp solution associated with the simulations and the reference solution presented above, the linearity of the differential operators and identical tangential $\bm{b}$-field boundary datum imply that
\begin{equation}
  \mu_0\partial_t r=\cdiv{\big(\kappa^{-1}\nabla r\big)},\qquad
  r=0\ \text{on the perimeter},\qquad r(x,y,0)=-g(x,y).
  \label{eq:square-startup-operator}
\end{equation}

Applying the relative error of Eq.~\eqref{eq:square-error-measures} to
the continuum current and using Eq.~\eqref{eq:square-reference-fields}
gives
\[
  \bar{e}^2(t)
  =\frac{\int_{\mathscr D}\kappa^{-1}|\nabla r|^2\,dv}
        {\int_{\mathscr D}\kappa^{-1}|\nabla g|^2\,dv}
  =\frac{\int\kappa^{-1}|\nabla r|^2\,da}
        {\int\kappa^{-1}|\nabla g|^2\,da}.
\]
The second equality follows because $r$, $g$, and $\kappa$ are
independent of $z$, so the thickness contributions cancel. Since $g$ is
independent of time, differentiation gives
\[
  \frac{\partial}{\partial t}\bar{e}^2(t)
  =2\bar{e}(t)[\bar{e}(t)]_{,t}
  =\frac{2\int\kappa^{-1}\nabla r\cdot\nabla r_{,t}\,da}
        {\int\kappa^{-1}|\nabla g|^2\,da}.
\]
Dividing both sides by $-2\bar{e}^2(t)$, applying the 2D-divergence theorem, and noting Eq.~\eqref{eq:square-startup-operator} then yields
\begin{equation}
\begin{aligned}
  -\frac{[\bar{e}(t)]_{,t}}{\bar{e}(t)}
  &=-\frac{\int\kappa^{-1}\nabla r\cdot\nabla r_{,t}\,da}
           {\int\kappa^{-1}|\nabla r|^2\,da}=\frac{\int\left[
       r_{,t}\,\cdiv{\big(\kappa^{-1}\nabla r\big)}
       -\cdiv{\big(r_{,t}\kappa^{-1}\nabla r\big)}
     \right]\,da}
     {\int\kappa^{-1}|\nabla r|^2\,da}=\frac{\mu_0\int(r_{,t})^2\,da}
          {\int\kappa^{-1}|\nabla r|^2\,da}.
\end{aligned}
\label{eq:square-diffusion-bound}
\end{equation}

To establish a bound on $\bar{e}(t)$ from above, it suffices to bound $-[\bar{e}(t)]_{,t}/\bar{e}(t)$ in Eq.~\eqref{eq:square-diffusion-bound} from below. First, note that Eq.~\eqref{eq:square-startup-operator} implies
\[
  \mu_0\int r\,r_{,t}\,da
  =\int r\,\cdiv{\big(\kappa^{-1}\nabla r\big)}\,da
  =-\int\kappa^{-1}|\nabla r|^2\,da .
\]
By Cauchy--Schwarz and the identity just obtained,
\[
  \left(\int\kappa^{-1}|\nabla r|^2\,da\right)^2
  =\mu_0^2\left(\int r\,r_{,t}\,da\right)^2
  \leq\mu_0^2\left(\int r^2\,da\right)\left(\int(r_{,t})^2\,da\right).
\]
Dividing both sides by $\mu_0\left(\int r^2\,da\right)\left(\int\kappa^{-1}|\nabla r|^2\,da\right)$ gives
\[
  \frac{\int\kappa^{-1}|\nabla r|^2\,da}{\mu_0\int r^2\,da}
  \leq\frac{\mu_0\int(r_{,t})^2\,da}{\int\kappa^{-1}|\nabla r|^2\,da}
  =-\frac{[\bar{e}(t)]_{,t}}{\bar{e}(t)},
\]
where the equality follows from Eq.~\eqref{eq:square-diffusion-bound}.

The boundary data allow us to bound the quotient on the left by a known constant. Indeed, since $r$ vanishes at both ends of each horizontal line through the square, the one-dimensional Poincar\'e inequality gives
\[
  \int r^2\,da\leq\frac{(1.5\times10^{-3})^2}{\pi^2}\int(r_{,x})^2\,da .
\]
Applying the same inequality along the vertical lines gives an analogous bound involving $r_{,y}$. Adding these two bounds, dividing by two, and using $\kappa\leq\kappa_a$ gives
\[
  \int r^2\,da
  \leq\frac{(1.5)^2\times10^{-6}}{2\pi^2}\int|\nabla r|^2\,da
  \leq\frac{1.125\times10^{-6}\kappa_a}{\pi^2}\int\kappa^{-1}|\nabla r|^2\,da .
\]
Rearranging yields the desired bound:
\begin{equation}
  \frac{10^6\pi^2}{1.125\mu_0\kappa_a}
  \leq\frac{\int\kappa^{-1}|\nabla r|^2\,da}{\mu_0\int r^2\,da}
  \leq-\frac{[\bar{e}(t)]_{,t}}{\bar{e}(t)} .
  \label{eq:square-transient-modes}
\end{equation}

The current starts from zero, so $\bar{e}(0)=1$. Integrating the right-hand side of Eq.~\eqref{eq:square-transient-modes} from zero to $t$ therefore gives $-\log \bar{e}(t)$, and hence
\begin{equation}
  \bar{e}(t)\leq\exp\!\left[-\frac{10^6\pi^2}{1.125\mu_0\kappa_a}t\right].
  \label{eq:square-transient-time}
\end{equation}
The time constant in this bound is $3.6354\,\mu\mathrm{s}$. Thus, by $25.2\,\mu\mathrm{s}$ the continuum current is within $0.1\%$ of the steady-ramp current in the norm of Eq.~\eqref{eq:square-current-norm}.

\end{appendices}

\end{document}